\font\fifteenrm=cmr10 scaled\magstep2 
\font\fifteeni=cmmi10 scaled\magstep2
\font\fifteensy=cmsy10 scaled\magstep2
\font\fifteenbf=cmbx10 scaled\magstep2
\font\fifteentt=cmtt10 scaled\magstep2
\font\fifteenit=cmti10 scaled\magstep2
\font\fifteensl=cmsl10 scaled\magstep2
\font\fifteenam=msam10 scaled\magstep2
\font\fifteenbm=msbm10 scaled\magstep2
\font\fifteenex=cmex10 scaled\magstep2
\font\fifteensc=cmcsc10 scaled\magstep2 
\font\twelverm=cmr10 at 12pt
\font\twelvei=cmmi10 at 12pt
\font\twelvesy=cmsy10 at 12pt
\font\twelvebf=cmbx10 at 12pt
\font\twelvett=cmtt10 at 12pt
\font\twelveit=cmti10 at 12pt
\font\twelvesl=cmsl10 at 12pt
\font\twelveam=msam10 at 12pt
\font\twelvebm=msbm10 at 12pt
\font\twelveex=cmex10 at 12pt
\font\twelvesc=cmcsc10 at 12pt
\font\elevenrm=cmr10 scaled\magstephalf 
\font\eleveni=cmmi10 scaled\magstephalf
\font\elevensy=cmsy10 scaled\magstephalf
\font\elevenbf=cmbx10 scaled\magstephalf
\font\eleventt=cmtt10 scaled\magstephalf
\font\elevenit=cmti10 scaled\magstephalf
\font\elevensl=cmsl10 scaled\magstephalf
\font\elevenam=msam10 scaled\magstephalf
\font\elevenbm=msbm10 scaled\magstephalf
\font\elevenex=cmex10 scaled\magstephalf
\font\elevensc=cmcsc10 scaled\magstephalf
\font\tenrm=cmr10
\font\teni=cmmi10
\font\tensy=cmsy10
\font\tenbf=cmbx10
\font\tentt=cmtt10
\font\tenit=cmti10
\font\tensl=cmsl10
\font\tenam=msam10
\font\tenbm=msbm10
\font\tenex=cmex10
\font\tensc=cmcsc10
\font\ninerm=cmr9
\font\ninei=cmmi9
\font\ninesy=cmsy9
\font\ninebf=cmbx9
\font\ninett=cmtt9
\font\nineit=cmti9
\font\ninesl=cmsl9
\font\nineam=msam9
\font\ninebm=msbm9
\font\nineex=cmex9
\font\ninesc=cmcsc9
\font\eightrm=cmr8
\font\eighti=cmmi8
\font\eightsy=cmsy8
\font\eightbf=cmbx8
\font\eighttt=cmtt8
\font\eightit=cmti8
\font\eightsl=cmsl8
\font\eightam=msam8
\font\eightbm=msbm8
\font\eightex=cmex8
\font\eightsc=cmcsc8
\font\sevenrm=cmr7
\font\seveni=cmmi7
\font\sevensy=cmsy7
\font\sevenbf=cmbx7

\font\sevenam=msam7
\font\sevenbm=msbm7

\font\sixrm=cmr6
\font\sixi=cmmi6
\font\sixsy=cmsy6

\font\sixam=msam6
\font\sixbm=msbm6

\font\fiverm=cmr5
\font\fivei=cmmi5
\font\fivesy=cmsy5

\font\fiveam=msam5
\font\fivebm=msbm5

\font\fourrm=cmr5 at 4pt
\font\fouri=cmmi5 at 4pt
\font\foursy=cmsy5 at 4pt

\font\fouram=msam5 at 4pt
\font\fourbm=msbm5 at 4pt

\skewchar\twelvei='177 \skewchar\eleveni='177\skewchar\teni='177
\skewchar\ninei='177 \skewchar\eighti='177\skewchar\seveni='177 
\skewchar\sixi='177 \skewchar\fivei='177 \skewchar\fouri='177
\skewchar\twelvesy='60 \skewchar\elevensy='60 \skewchar\tensy='60
\skewchar\ninesy='60 \skewchar\eightsy='60 \skewchar\sevensy='60 
\skewchar\sixsy='60 \skewchar\fivesy='60 \skewchar\foursy='60
\newfam\itfam
\newfam\slfam
\newfam\bffam
\newfam\ttfam
\newfam\scfam
\newfam\amfam
\newfam\bmfam
\def\eightbig#1{{\hbox{$\left#1\vbox to 6.5pt{}\voidright $}}}
\def\eightBig#1{{\hbox{$\left#1\vbox to 7.5pt{}\voidright $}}}
\def\eightbigg#1{{\hbox{$\left#1\vbox to 10pt{}\voidright $}}}
\def\eightBigg#1{{\hbox{$\left#1\vbox to 13pt{}\voidright $}}}
\def\ninebig#1{{\hbox{$\left#1\vbox to 7.5pt{}\voidright $}}}
\def\nineBig#1{{\hbox{$\left#1\vbox to 8.5pt{}\voidright $}}}
\def\ninebigg#1{{\hbox{$\left#1\vbox to 11.5pt{}\voidright $}}}
\def\nineBigg#1{{\hbox{$\left#1\vbox to 14.5pt{}\voidright $}}}
\def\tenbig#1{{\hbox{$\left#1\vbox to 8.5pt{}\voidright $}}}
\def\tenBig#1{{\hbox{$\left#1\vbox to 9.5pt{}\voidright $}}}
\def\tenbigg#1{{\hbox{$\left#1\vbox to 12.5pt{}\voidright $}}}
\def\tenBigg#1{{\hbox{$\left#1\vbox to 16pt{}\voidright $}}}
\def\elevenbig#1{{\hbox{$\left#1\vbox to 9pt{}\voidright $}}}
\def\elevenBig#1{{\hbox{$\left#1\vbox to 10.5pt{}\voidright $}}}
\def\elevenbigg#1{{\hbox{$\left#1\vbox to 14pt{}\voidright $}}}
\def\elevenBigg#1{{\hbox{$\left#1\vbox to 17.5pt{}\voidright $}}}
\def\twelvebig#1{{\hbox{$\left#1\vbox to 10pt{}\voidright $}}}
\def\twelveBig#1{{\hbox{$\left#1\vbox to 11pt{}\voidright $}}}
\def\twelvebigg#1{{\hbox{$\left#1\vbox to 15pt{}\voidright $}}}
\def\twelveBigg#1{{\hbox{$\left#1\vbox to 19pt{}\voidright $}}}
\def\fifteenbig#1{{\hbox{$\left#1\vbox to 12pt{}\voidright $}}}
\def\fifteenBig#1{{\hbox{$\left#1\vbox to 13.5pt{}\voidright $}}}
\def\fifteenbigg#1{{\hbox{$\left#1\vbox to 18pt{}\voidright $}}}
\def\fifteenBigg#1{{\hbox{$\left#1\vbox to 23pt{}\voidright $}}}
\def\voidright{\right.\nulldelimiterspace=0pt \mathsurround=0pt }
\def\fifteenpoint{
  \textfont0=\fifteenrm \scriptfont0=\twelverm \scriptscriptfont0=\tenrm
  \def\rm{\fam0 \fifteenrm}%
  \textfont1=\fifteeni \scriptfont1=\twelvei \scriptscriptfont1=\teni
  \textfont2=\fifteensy \scriptfont2=\twelvesy \scriptscriptfont2=\tensy
  \textfont3=\fifteenex \scriptfont3=\fifteenex \scriptscriptfont3=\fifteenex
  \def\it{\fam\itfam\fifteenit}\textfont\itfam=\fifteenit
  \def\sl{\fam\slfam\fifteensl}\textfont\slfam=\fifteensl
  \def\bf{\fam\bffam\fifteenbf}\textfont\bffam=\fifteenbf 
    \scriptfont\bffam=\twelvebf\scriptscriptfont\bffam=\tenbf
  \def\tt{\fam\ttfam\fifteentt}\textfont\ttfam=\fifteentt
  \def\sc{\fam\scfam\fifteensc}\textfont\scfam=\fifteensc
  \def\am{\fam\amfam\fifteenam}\textfont\amfam=\fifteenam
    \scriptfont\amfam=\twelveam\scriptscriptfont\amfam=\tenam
  \def\bm{\fam\bmfam\fifteenbm}\textfont\bmfam=\fifteenbm
    \scriptfont\bmfam=\twelvebm\scriptscriptfont\bmfam=\tenbm
  \baselineskip=21pt \rm
  \let\big=\fifteenbig\let\Big=\fifteenBig\let\bigg=\fifteenbigg
  \let\Bigg=\fifteenBigg}
\def\twelvepoint{
  \textfont0=\twelverm \scriptfont0=\ninerm \scriptscriptfont0=\sevenrm
  \def\rm{\fam0 \twelverm}%
  \textfont1=\twelvei \scriptfont1=\ninei \scriptscriptfont1=\seveni
  \textfont2=\twelvesy \scriptfont2=\ninesy \scriptscriptfont2=\sevensy
  \textfont3=\twelveex \scriptfont3=\twelveex \scriptscriptfont3=\twelveex
  \def\it{\fam\itfam\twelveit}\textfont\itfam=\twelveit
  \def\sl{\fam\slfam\twelvesl}\textfont\slfam=\twelvesl
  \def\bf{\fam\bffam\twelvebf}\textfont\bffam=\twelvebf 
    \scriptfont\bffam=\ninebf\scriptscriptfont\bffam=\sevenbf
  \def\tt{\fam\ttfam\twelvett}\textfont\ttfam=\twelvett
  \def\sc{\fam\scfam\twelvesc}\textfont\scfam=\twelvesc
  \def\am{\fam\amfam\twelveam}\textfont\amfam=\twelveam
    \scriptfont\amfam=\nineam\scriptscriptfont\amfam=\sevenam
  \def\bm{\fam\bmfam\twelvebm}\textfont\bmfam=\twelvebm
    \scriptfont\bmfam=\ninebm\scriptscriptfont\bmfam=\sevenbm
  \baselineskip=17.8pt \rm 
  \def\looselineskip{\baselineskip=18.5pt plus 1.8pt}%
  \def\tightlineskip{\baselineskip=16.5pt}%
  \def\verytightlineskip{\baselineskip=15pt}%
  \let\big=\twelvebig\let\Big=\twelveBig\let\bigg=\twelvebigg
  \let\Bigg=\twelveBigg  }
\def\elevenpoint{
  \textfont0=\elevenrm \scriptfont0=\ninerm \scriptscriptfont0=\sixrm
  \def\rm{\fam0 \elevenrm}%
  \textfont1=\eleveni \scriptfont1=\ninei \scriptscriptfont1=\sixi
  \textfont2=\elevensy \scriptfont2=\ninesy \scriptfont2=\sixsy 
  \textfont3=\elevenex \scriptfont3=\elevenex \scriptfont3=\elevenex
  \def\it{\fam\itfam\elevenit}\textfont\itfam=\elevenit
  \def\sl{\fam\slfam\elevensl}\textfont\slfam=\elevensl
  \def\bf{\fam\bffam\elevenbf}\textfont\bffam=\elevenbf
  \def\tt{\fam\ttfam\eleventt}\textfont\ttfam=\eleventt
  \def\sc{\fam\scfam\elevensc}\textfont\scfam=\elevensc
  \def\am{\fam\amfam\elevenam}\textfont\amfam=\elevenam
    \scriptfont\amfam=\nineam\scriptscriptfont\amfam=\sixam
  \def\bm{\fam\bmfam\elevenbm}\textfont\bmfam=\elevenbm
    \scriptfont\bmfam=\ninebm\scriptscriptfont\bmfam=\sixbm
  \baselineskip=15.1pt \rm
  \def\looselineskip{\baselineskip=16pt plus 1.5pt}%
  \def\tightlineskip{\baselineskip=14pt}%
  \def\verytightlineskip{\baselineskip=13pt}%
  \let\big=\elevenbig\let\Big=\elevenBig\let\bigg=\elevenbigg
  \let\Bigg=\elevenBigg  }
\def\tenpoint{
  \textfont0=\tenrm \scriptfont0=\eightrm \scriptscriptfont0=\fiverm
  \def\rm{\fam0 \tenrm}%
  \textfont1=\teni \scriptfont1=\eighti \scriptscriptfont1=\fivei
  \textfont2=\tensy \scriptfont2=\eightsy \scriptfont2=\fivesy 
  \textfont3=\tenex \scriptfont3=\tenex \scriptfont3=\tenex
  \def\it{\fam\itfam\tenit}\textfont\itfam=\tenit
  \def\sl{\fam\slfam\tensl}\textfont\slfam=\tensl
  \def\bf{\fam\bffam\tenbf}\textfont\bffam=\tenbf
  \def\tt{\fam\ttfam\tentt}\textfont\ttfam=\tentt
  \def\sc{\fam\scfam\tensc}\textfont\scfam=\tensc
  \def\am{\fam\amfam\tenam}\textfont\amfam=\tenam
    \scriptfont\amfam=\eightam \scriptscriptfont\amfam=\fiveam
  \def\bm{\fam\bmfam\tenbm}\textfont\bmfam=\tenbm
    \scriptfont\bmfam=\eightbm \scriptscriptfont\bmfam=\fivebm
  \baselineskip=14pt \rm
  \def\looselineskip{\baselineskip=14.8pt plus1.5pt}
  \def\tightlineskip{\baselineskip=13.6pt}%
  \def\verytightlineskip{\baselineskip=13pt}%
  \let\big=\tenbig\let\Big=\tenBig\let\bigg=\tenbigg\let\Bigg=\tenBigg  }
\def\ninepoint{
  \textfont0=\ninerm \scriptfont0=\sevenrm \scriptscriptfont0=\fourrm
  \def\rm{\fam0 \ninerm}%
  \textfont1=\ninei \scriptfont1=\seveni \scriptscriptfont1=\fouri
  \textfont2=\ninesy \scriptfont2=\sevensy \scriptfont2=\foursy 
  \textfont3=\nineex \scriptfont3=\nineex \scriptfont3=\nineex
  \def\it{\fam\itfam\nineit}\textfont\itfam=\nineit
  \def\sl{\fam\slfam\ninesl}\textfont\slfam=\ninesl
  \def\bf{\fam\bffam\ninebf}\textfont\bffam=\ninebf
  \def\tt{\fam\ttfam\ninett}\textfont\ttfam=\ninett
  \def\sc{\fam\scfam\ninesc}\textfont\scfam=\ninesc
  \def\am{\fam\amfam\nineam}\textfont\amfam=\nineam
    \scriptfont\amfam=\nineam\scriptscriptfont\amfam=\fouram
  \def\bm{\fam\bmfam\ninebm}\textfont\bmfam=\ninebm
    \scriptfont\bmfam=\ninebm\scriptscriptfont\bmfam=\fourbm
  \baselineskip=12.6pt \rm
  \let\big=\ninebig\let\Big=\nineBig\let\bigg=\ninebigg
  \let\Bigg=\nineBigg  }
\def\eightpoint{
  \textfont0=\eightrm \scriptfont0=\fiverm \scriptscriptfont0=\fourrm
  \def\rm{\fam0 \eightrm}%
  \textfont1=\eighti \scriptfont1=\fivei \scriptscriptfont1=\fouri
  \textfont2=\eightsy \scriptfont2=\fivesy \scriptfont2=\foursy 
  \textfont3=\eightex \scriptfont3=\eightex \scriptfont3=\eightex
  \def\it{\fam\itfam\eightit}\textfont\itfam=\eightit
  \def\sl{\fam\slfam\eightsl}\textfont\slfam=\eightsl
  \def\bf{\fam\bffam\eightbf}\textfont\bffam=\eightbf
  \def\tt{\fam\ttfam\eighttt}\textfont\ttfam=\eighttt
  \def\sc{\fam\scfam\eightsc}\textfont\scfam=\eightsc
  \def\am{\fam\amfam\eightam}\textfont\amfam=\eightam
    \scriptfont\amfam=\eightam\scriptscriptfont\amfam=\fouram
  \def\bm{\fam\bmfam\eightbm}\textfont\bmfam=\eightbm
    \scriptfont\bmfam=\eightbm\scriptscriptfont\bmfam=\fourbm
  \baselineskip=11.2pt \rm
  \let\big=\eightbig\let\Big=\eightBig\let\bigg=\eightbigg
  \let\Bigg=\eightBigg  }

\twelvepoint
\nopagenumbers
\hsize=6in\vsize=8.8in

\parskip=1pt plus 1pt

\newif\ifSpecialhead\Specialheadfalse
\newbox\specialheadbox

\def\specialhead #1\par{\Specialheadtrue\setbox\specialheadbox=\hbox{#1}}
\headline={{\ifSpecialhead\box\specialheadbox\global\Specialheadfalse\else
     \ifnum\pageno<0{\hfill\quad{\twelvebf\folio}}%
     \else\ifnum\pageno<2\hfill
     \else\hfill\twelvepoint\sc\firstmark\quad{\twelvebf\folio}\fi\fi\fi}}

\def\title#1\par{\bigskip{\def\cr{\par\center}\center\fifteenbf #1\par}\medskip}
\def\subtitle#1\par{\centerline{\fifteenrm #1}\medskip}
\def\author#1\par{\medskip{\def\cr{\par\center\twelvesc}\fifteensc\center#1\par}}
\def\center#1\par{\hfil #1\hfil\par}
\def\abstract.#1\par{\message{Abstract.}%
                    \medskip{\narrower\narrower\tenpoint\tightlineskip
                        \noindent{\bf Abstract.}#1\par}\medskip\noindent}
\def\bigabstract.#1\par{\message{Abstract.}%
                         \medskip{\narrower\narrower\tightlineskip
                         \noindent{\bf Abstract. }#1\par}\medskip\noindent}
\def\acknowledgement#1\par{\footnote{}{#1}}
\def\sectionskip{\Goodbreak\vskip 25pt plus 15pt minus 5pt}
\def\secnumber{\ifquiet
               \else\ifNoSections
                    \else\sectionsymbol\the\secno\quad\fi\fi}
\def\section#1\par{ \NoSectionsfalse\par\sectionskip\proofdepth=0\claimno=0
 \ifquiet\else\advance\secno by1\fi\toks0={#1}
 \immediate\write16{\ifquiet\else Section \the\secno\space\fi
                    \the\toks0}%
 \mark{\secnumber #1}%
 {\fifteenpoint\bf\noindent\secnumber #1}\nobreak\bigskip\quietoff
 \nobreak\noindent}

\def\QUIET{\QUIETtrue\quiettrue}

\def\quietoff{\ifQUIET\else\quietfalse\fi}
\newif\ifquiet
\newif\ifQUIET
\newif\ifNoSections
\newcount\claimtype
\newcount\secno
\newcount\claimno
\newcount\subclaimno
\newcount\subsubclaimno
\newcount\subsubsubclaimno
\newcount\proofdepth
\def\subclaimnumber{\ifquiet\else\ifcase\subclaimno\or A\or B\or C\or D\or E\or
     F\or G\or H\or I\or J\or K\or L\or M\or N\or O\or P\fi\fi}
\def\subsubclaimnumber{\ifquiet\else\ifcase\subsubclaimno\or i\or ii\or iii\or 
   iv\or v\or vi\or vii\or viii\or ix\or x\or xi\or xii\or xiii\or xiv\fi\fi}
\def\subsubsubclaimnumber{\ifquiet\else\ifcase\subsubsubclaimno\or a\or b\or 
   c\or d\or e\or f\or g\or \or h\or i\or j\or k\or l\or m\or n\or o\fi\fi}
\def\claimtag{\ifquiet\else
  \ifNoSections
    \ifcase\proofdepth\the\claimno%
    \or\the\claimno.\subclaimnumber
    \or\the\claimno.\subclaimnumber.\subsubclaimnumber
    \or\the\claimno.\subclaimnumber.\subsubclaimnumber
                                                .\subsubsubclaimnumber\fi
  \else
    \ifcase\proofdepth\the\secno.\the\claimno
    \or\the\secno.\the\claimno.\subclaimnumber
    \or\the\secno.\the\claimno.\subclaimnumber.\subsubclaimnumber
    \or\the\secno.\the\claimno.\subclaimnumber.\subsubclaimnumber
                                                .\subsubsubclaimnumber\fi\fi\fi}
\secno=0\claimno=0\proofdepth=0\subclaimno=0\subsubclaimno=0\subsubsubclaimno=0
\NoSectionstrue
\newbox\qedbox
\def\claimname{\ifcase\claimtype Theorem\or Lemma\or Claim\or Corollary\or
               Question\or Definition\or Remark\or Conjecture\fi}
\def\preclaimskip{\removelastskip
    \ifcase\claimtype\goodbreak\vskip 8pt plus 4pt minus 2pt
                  \or\goodbreak\vskip 6pt plus 4pt minus 1pt
                  \or\goodbreak\vskip 5pt plus 4pt minus 1pt
                  \or\goodbreak\vskip 8pt plus 4pt minus 2pt
                  \or\vskip 7pt plus 4pt minus 2pt
                  \or\vskip 7pt plus 4pt minus 2pt
                  \or\vskip 7pt plus 4pt minus 2pt
                  \or\goodbreak\vskip 8pt plus 4pt minus 2pt\fi}
\def\postclaimskip{\ifcase\claimtype         \vskip 4pt plus 2pt minus 2pt
                                          \or\vskip 3pt plus 2pt minus 2pt
                                          \or\vskip 2pt plus 2pt minus 1pt
                                          \or\vskip 4pt plus 2pt minus 2pt
                                          \or\vskip 1pt plus 2pt 
                                          \or\vskip 4pt plus 4pt 
                                          \or\vskip 3pt plus 2pt
                                          \or\vskip 4pt plus 2pt minus 2pt\fi}
\def\claimfont{\ifcase\claimtype
                  \sl\or\sl\or\sl\or\sl\or\sl\or\rm\or\rm\or\sl\fi}
\def\advancetag{\ifcase\proofdepth\advance\claimno by1
                               \or\advance\subclaimno by1
                               \or\advance\subsubclaimno by1
                               \or\advance\subsubsubclaimno by1\fi}
\def\sayclaim#1.#2 #3\par{\ifquiet\else\advancetag\fi
    \preclaimskip\setbox1=\hbox{#1}\setbox2=\hbox{#2}%
    \toks0={#1 }
    \immediate\write16{\ifdim\wd1>0pt\the\toks0
                       \else\claimname\space\fi \claimtag.}%
    \vbox{\noindent
    {\bf\ifdim\wd1=0pt \claimname\else #1\fi\ifquiet.\else\ \claimtag{\ifNoSections.\fi}\fi}%
    \enspace{\ifdim\wd2>0pt\sc #2\enspace\fi}%
    {\claimfont #3\par}}\postclaimskip\quietoff}
\def\theorem{\claimtype=0\sayclaim}

\def\point#1. #2\par{\item{\rm #1.}#2\par}
\def\points#1\cr\par{\medskip\vbox{\let\cr=\point\point#1\par}\par}

\def\prooffont{}
\def\proofsize{}
\def\proofindent{}
\def\proofskip{\badbreak\ifcase\claimtype    \vskip 3pt plus 2pt minus 2pt
                                          \or\vskip 2pt plus 2pt minus 2pt
                                          \or\vskip 1pt plus 2pt minus 1pt
                                          \or\vskip 3pt plus 2pt minus 2pt
                                          \or\vskip 1pt plus 2pt 
                                          \or\vskip 2pt plus 4pt 
                                          \or\vskip 1pt plus 2pt
                                          \or\vskip 3pt plus 2pt minus 2pt\fi}

\def\Goodbreak{\vskip0pt plus.5in\penalty-1000\vskip0pt plus-.5in}
\def\goodbreak{\penalty-500}
\def\badbreak{\penalty500}
\def\Badbreak{\penalty1000}
\def\proof{\message{proof}\removelastskip\Badbreak\proofskip\begingroup
  \advance\proofdepth by1
  \setbox\qedbox=\hbox{\halmos\raise2pt\hbox{\fiverm\claimname}}%
  \prooffont\proofsize\proofindent\noindent{\bf Proof: }}
\def\proofof#1:{\message{proof}\removelastskip\Badbreak\proofskip\begingroup
  \advance\proofdepth by1
  \setbox\qedbox=\hbox{\halmos\raise2pt\hbox{\fiverm#1}}%
  \prooffont\proofsize\proofindent\noindent{\bf Proof of #1: }}
\def\cite[#1]{[{\tenrm{#1}}]\message{[#1]}}
\edef\ref#1{\expandafter\global\expandafter\edef#1{\noexpand\claimtag}}
\newwrite\notes
\openout\notes=\jobname.notes
\long\def\unexpandedwrite#1#2{\def\finwrite{\write#1}%
   {\aftergroup\finwrite\aftergroup{\sanitize#2\endsanity}}}
\def\sanitize{\futurelet\next\sanswitch}
\let\stoken=\space
\def\sanswitch{\ifx\next\endsanity
  \else\ifcat\noexpand\next\stoken\aftergroup\space\let\next=\eat
   \else\ifcat\noexpand\next\bgroup\aftergroup{\let\next=\eat
    \else\ifcat\noexpand\next\egroup\aftergroup}\let\next=\eat
     \else\let\next=\copytoken\fi\fi\fi\fi \next}
\def\eat{\afterassignment\sanitize \let\next= }
\long\def\copytoken#1{\ifcat\noexpand#1\relax\aftergroup\noexpand
  \else\ifcat\noexpand#1\noexpand~\aftergroup\noexpand\fi\fi
  \aftergroup#1\sanitize}
\def\endsanity\endsanity{}

\def\note#1#2{\hbox to2in{\strut#1\quad\dotfill\quad#2}}
\def\boxit#1{\setbox4=\hbox{\kern1pt#1\kern1pt}
  \hbox{\vrule\vbox{\hrule\kern1pt\box4\kern1pt\hrule}\vrule}}
\def\halmos{\hbox{\am\char'3}} 
\def\qed#1\par{\message{.                                }\setbox1=\hbox{#1}%
  \ifdim\wd1>0pt\setbox\qedbox=\hbox{\halmos\raise2pt\hbox{\fiverm #1}}\fi
  \kern5pt\lower 2pt\hbox{\box\qedbox}\proofskip\goodbreak\endgroup}

\def\sectionsymbol{\S}
\def\k{\kappa}
\def\g{\gamma}
\def\a{\alpha}
\def\b{\beta}

\def\l{\lambda}

\def\I1{\mathop{\hbox{\sc i}_1}}
\def\w{\omega}

\def\Z{{\mathchoice{\hbox{\bm Z}}{\hbox{\bm Z}}
         {\hbox{\tenbm Z}}{\hbox{\sevenbm Z}}}}
\def\N{{\mathchoice{\hbox{\bm N}}{\hbox{\bm N}}
         {\hbox{\tenbm N}}{\hbox{\sevenbm N}}}}

\def\card#1{\left|#1\right|}

\def\dirlim{\mathop{\rm dir\,lim}\nolimits}

\def\aut{\mathop{\hbox{\rm Aut}}\nolimits}

\def\unifto{\buildrel\lower 7pt\hbox{$\to$}\over\to}

\def\iso{\cong}

\def\plus{^{\scriptscriptstyle +}}

\def\in{\mathrel{\mathchoice{\raise 
1pt\hbox{$\scriptstyle\cal\char'62$}}
         {\raise 1pt\hbox{$\scriptstyle\cal\char'62$}}
         {\raise .5pt\hbox{$\scriptscriptstyle\cal\char'62$}}
         {\hbox{$\scriptscriptstyle\cal\char'62$}}}\penalty700{}}
\def\ni{\mathrel{\mathchoice{\raise 1pt\hbox{$\scriptstyle\cal\char'63$}}
                   {\raise 1pt\hbox{$\scriptstyle\cal\char'63$}}
                   {\raise .5pt\hbox{$\scriptscriptstyle\cal\char'63$}}
                   {\hbox{$\scriptscriptstyle\cal\char'63$}}}\penalty700}
\def\of{\mathrel{\mathchoice{\raise 1pt\hbox{$\scriptstyle\subseteq$}}
                   {\raise 1pt\hbox{$\scriptstyle\subseteq$}}
                   {\raise .5pt\hbox{$\scriptscriptstyle\subseteq$}}
                   {\hbox{$\scriptscriptstyle\subseteq$}}}}
\def\fo{\mathrel{\mathchoice{\raise 1pt\hbox{$\scriptstyle\supseteq$}}
                   {\raise 1pt\hbox{$\scriptstyle\supseteq$}}
                   {\raise .5pt\hbox{$\scriptscriptstyle\supseteq$}}
                   {\hbox{$\scriptscriptstyle\supseteq$}}}}
\def\notin{\mathrel{\mathchoice
  {\raise 1pt\hbox{\rlap{$\scriptstyle\;|$}$\scriptstyle\cal\char'62$}}
  {\raise 1pt\hbox{\rlap{$\scriptstyle\kern2pt 
          |$}$\scriptstyle\cal\char'62$}}
  {\raise .5pt\hbox{\rlap{$\scriptscriptstyle\, |$}$\scriptscriptstyle
      \cal\char'62$}}
  {\hbox{\rlap{$\scriptscriptstyle\, |$}$\scriptscriptstyle
     \cal\char'62$}}}%
  \penalty700}

\def\and{\mathrel{\kern1pt\&\kern1pt}}

\def\cross{\times}

\def\[#1]{\left[\vphantom{\bigm|}#1\right]}
\def\<#1>{\langle\,#1\,\rangle}

\def\st{\mid}
\def\seq<#1>{{\def\st{\mid\penalty650}\left<\,#1\,\right>}}

\def\set#1{\{\,#1\,\}}

\def\lttheta{{\raise 1pt\hbox{$\scriptstyle<$}\theta}}

\def\I1{\mathop{\hbox{\sc i}_1}}

\font\arrow=line10 scaled \magstep1
\def\makeline#1.{\hbox{\arrow\char#1}}
\def\makearrow#1.#2.{\hbox{\arrow\char#1\llap{\char#2}}}
\def\definelinesandarrows#1.#2.#3.#4.#5.{
   \expandafter\edef\csname#4line\endcsname{\makeline#1.}
   \expandafter\edef\csname#4arrow\endcsname{\makearrow#1.#2.}
   \expandafter\edef\csname#5line\endcsname{\makeline#1.}
   \expandafter\edef\csname#5arrow\endcsname{\makearrow#1.#3.}}
\definelinesandarrows 0.18.9.ne.sw.
\definelinesandarrows 1.21.11.nnne.sssw.
\definelinesandarrows 2.14.13.nnnne.ssssw.
\definelinesandarrows 3.23.15.nnnnne.sssssw.
\definelinesandarrows 4.23.15.nnnnnne.ssssssw.
\definelinesandarrows 10.30.29.nne.ssw.
\definelinesandarrows 16.49.41.neeeeee.swwwwww.
\definelinesandarrows 17.51.43.neeee.swwww.
\definelinesandarrows 19.55.47.nehuh.swhuh.
\definelinesandarrows 24.58.41.neeeeeee.swwwwwww.
\definelinesandarrows 26.62.9.neee.swww.
\definelinesandarrows 33.49.25.neeeee.swwwww.
\definelinesandarrows 35.62.61.nee.sww.
\definelinesandarrows 64.82.73.se.nw.
\definelinesandarrows 65.85.75.ssse.nnnw.
\definelinesandarrows 66.78.77.sssse.nnnnw.
\definelinesandarrows 67.87.79.ssssse.nnnnnw.
\definelinesandarrows 68.87.79.sssssse.nnnnnnw.
\definelinesandarrows 74.94.93.sse.nnw.
\definelinesandarrows 80.113.105.seeeeee.nwwwwww.
\definelinesandarrows 81.115.107.seeee.nwwww.
\definelinesandarrows 99.126.125.see.nww.
\def\sejoin#1#2{\setbox1=\hbox{#1}\setbox2=\hbox{#2}%
  \hbox{\vbox{\hbox{\copy1\kern\wd2}\nointerlineskip
              \hbox{\kern\wd1\box2}}}}
\def\nejoin#1#2{\setbox1=\hbox{#1}\setbox2=\hbox{#2}%
  \hbox{\vbox{\hbox{\kern\wd1\copy2}\nointerlineskip\hbox{\copy1\kern\wd2}}}}
\newdimen\hnudge
\newdimen\vnudge
\newdimen\hnudgedefault
\newdimen\vnudgedefault

\def\SEdefaultnudge{\hnudge=-16pt\vnudge=20pt}
\def\Edefaultnudge{\hnudge=-25pt\vnudge=6pt}

\def\longEdefaultnudge{\hnudge=-5pt\vnudge=6pt}
\def\nudgeright#1pt{\advance\hnudge by#1pt}
\def\nudgeleft#1pt{\advance\hnudge by-#1pt}
\def\nudgeup#1pt{\advance\vnudge by#1pt}
\def\nudgedown#1pt{\advance\vnudge by-#1pt}
\def\label#1{\smash{\llap{\kern\hnudge
                   \raise\vnudge\rlap{$\scriptstyle#1$}\hfill}}}

\def\SEarrow{\SEdefaultnudge
             \sejoin\seeline{\sejoin\seeline{\sejoin\seeline\seearrow}}}

\def\Earrow{\Edefaultnudge\setbox1=\hbox{\SEarrow}
 \hbox{\raise 2pt\hbox{\vrule height-.4pt depth.8ptwidth\wd1\kern2pt
       \llap{\arrow\char'55}}}}
\def\longEarrow{\longEdefaultnudge\setbox1=\hbox{\SEarrow}
      \rlap{\hskip-1.25\wd1\raise 2pt
            \hbox{\vrule height-.4pt depth.8ptwidth2.5\wd1\kern2pt
            \llap{\arrow\char'55}}}}

\QUIET\twelvepoint\tightlineskip
\center {\tenrm [to appear in the Proceedings of the American 
    Mathematical Society]}

\title Every Group Has A Terminating\cr
Transfinite Automorphism Tower\cr

\author Joel David Hamkins\footnote\dag{{\rm My research has been
supported in part by a grant from the PSC-CUNY Research Foundation. I
would like to thank both Daniel Seabold and Daniel Velleman for
pointing out a simplification in my proof.}}\cr
	City University of New York\cr

\abstract. Iteratively taking the automorphism group of any group
leads, transfinitely, to a fixed point.

\noindent The automorphism tower of a group is obtained by computing
its automorphism group, the automorphism group of {\it that} group, and
so on, iterating transfinitely. Each group maps canonically into the
next using inner automorphisms, and so at limit stages one can take a
direct limit and continue the iteration.  $$G\to\aut(G)\to
\aut(\aut(G))\to\cdots\to G_\omega\to G_{\omega+1}\to \cdots\to
G_\a\to\cdots$$ The tower is said to terminate if a fixed point is
reached, that is, if a group is reached which is isomorphic to its
automorphism group by the natural map. This occurs if a {\it complete}
group is reached, one which is centerless and has only inner
automorphisms.

In the special case that the initial group is centerless, matters
simplify considerably: in this case {\it all} the groups appearing in
the tower are centerless (see Hulse [1970]), and, consequently, all the
natural maps are injective. The tower can therefore be viewed as
building upwards to larger and larger groups; the question is whether
this building process ever stops\footnote\ddag{In Hulse [1970], Rae and
Roseblade [1970], and Thomas [1985], the tower is only defined in this
special case; but the definition I give here works perfectly well
whether or not the group $G$ is centerless. Of course, when there is a
center, one has homomorphisms rather than embeddings.}.  Wielandt
[1939] proved the classical result that the automorphism tower of any
centerless finite group terminates in finitely many steps.  Rae and
Roseblade [1970] proved that the automorphism tower of any centerless
\v Cernikov group terminates in finitely many steps. Hulse [1970]
proved that the the automorphism tower of any centerless polycyclic
group terminates in countably many steps (but not necessarily after
just $\omega$ many steps). Solving the problem for centerless groups,
Simon Thomas [1985] proved that the automorphism tower of any
centerless group eventually terminates. In fact, the automorphism tower
of a centerless group $G$ terminates in fewer than
$(2^{\card{G}})\plus$ many steps.  In the general case, however, the
question remained open whether every group has a terminating
automorphism tower.  This is settled by the following theorem.

\theorem Main Theorem. Every group has a terminating automorphism tower.

\proof Suppose $G$ is a group. The following transfinite 
recursion defines the automorphism tower of $G$:
$$\eqalign{G_0=&\,G\cr
           G_{\a+1}=&\aut(G_\a), \quad
           \hbox{where $\pi_{\a,\a+1}:G_\a\to G_{\a+1}$ is the natural map,}\cr
           G_\l=&\dirlim_{\a<\l}G_\a, \quad
            \hbox{if $\l$ is a limit ordinal.}\cr}$$
When $\a<\b$ one obtains the map $\pi_{\a,\b}:G_\a\to G_\b$ by
composing the canonical maps at each step, and these are the maps used
to compute the direct limit at limit stages. Thus, when $\l$ is a limit
ordinal, every element of $G_\l$ is of the form $\pi_{\a,\l}(g)$ for
some $\a<\l$ and some $g\in G_\a$.

Since Simon Thomas [1985] has proved that every centerless group has a
terminating automorphism tower, it suffices to show that there is an
ordinal $\g$ such that $G_\g$ has a trivial center. For each ordinal
$\a$, let $H_\a=\set{g\in G_\a\st \exists\b\,\pi_{\a,\b}(g)=1}$.  For
every $g\in H_\a$ there is some least $\b_g>\a$ such that
$\pi_{\a,\b_g}(g)=1$.  Let $f(\a)=\sup_{g\in H_\a}\b_g$. It is easy to
check that if $\a<\b$ then $\a<f(\a)\leq f(\b)$. Iterating the
function, define $\g_0=0$ and $\g_{n+1}=f(\g_n)$. This produces a
strictly increasing $\w$-sequence of ordinals whose supremum
$\g=\sup\set{\g_n\st n\in\w}$ is a limit ordinal which is closed under
$f$. That is, $f(\a)<\g$ for every $\a<\g$.  I claim that $G_\g$ has a
trivial center. To see this, suppose $g$ is in the center of $G_\g$.
Thus, $\pi_{\g,\g+1}(g)=1$. Moreover, since $\g$ is a limit ordinal,
there is $\a<\g$ and $h\in G_\a$ such that $g=\pi_{\a,\g}(h)$.
Combining these facts, observe that
$$\pi_{\a,\g+1}(h)=\pi_{\g,\g+1}(\pi_{\a,\g}(h))=\pi_{\g,\g+1}(g)=1.$$
Consequently, $\pi_{\a,f(\a)}(h)=1$.  Since $f(\a)<\g$, it follows that
$$g=\pi_{\a,\g}(h)=\pi_{f(\a),\g}(\pi_{\a,f(\a)}(h))=\pi_{f(\a),\g}(1)=1,$$
as desired.\qed

The proof does not reveal exactly how long the automorphism tower takes
to stabilize, since it is not clear how large $f(\a)$ can be.
Nevertheless, something more can be said. Certainly the automorphism
tower of $G$ terminates well before the next inaccessible cardinal
above $\card{G}$.  More generally, if $\k>\w$ is regular and
$\card{G_\a}<\k$ whenever $\a<\k$, then I claim the centerless groups
will appear before $\k$. To see this, let $H_\a=\set{g\in G_\a\st
\exists\b_g{<}\k\,\, \pi_{\a,\b_g}(g)=1}$ and define $f:\k\to\k$ by
$f(\a)=\sup_{g\in H_\a}\b_g$; if $\g<\k$ is closed under $f$, then it
follows as in the main theorem that $G_\g$ has no center, as desired.
In this case the tower therefore terminates in fewer than
$(2^{\card{G_\g}})\plus$ many additional steps. If it happens that
$\card{G_\g}\plus<\k$, one can adapt Thomas' \cite[1996] argument using
Fodor's lemma to prove that the tower terminates actually in fewer than
$\k$ many steps. The point is that one can find a bound on the height
of the tower by bounding the rate of growth of the groups in the
tower.

Thomas \cite[1985] provides his explicit bound in the case of
centerless $G$ in precisely this way. He proves that if $G$ is
centerless and $\k=(2^{\card{G}})\plus$, then $\card{G_\a}<\k$ for all
$\a<\k$. The analogous result, unfortunately, does not hold for groups
with nontrivial centers. This is illustrated by the following example,
provided by the anonymous referee of this paper:

\theorem Example. There exists a countable group $G$ such that 
$\card{\aut G}=2^\w$ and\break $\card{\aut(\aut G)}=2^{2^\w}$.

\proof For each prime $p$, let $\Z[1/p]=\set{m/{p^n}\st m\in\Z,n\in\N}$
be the additive group of $p$-adic rationals and let $G=\oplus_p\Z[1/p]$
be the direct sum of these groups. An element $g\in G$ is divisible by $p^n$
for all $n\in\N$ iff $g\in\Z[1/p]$. Hence, if $\pi\in\aut G$, then
$\pi[\Z[1/p]]=\Z[1/p]$ for each prime $p$. It is easy to see that any
automorphism of $\Z[1/p]$ is simply multiplication by an element $u\in
U_p=\set{\pm p^n\st n\in\Z}$, the group of multiplicative units of the
ring of $p$-adic rationals. Thus, $\aut G\iso\Pi_p U_p$; and so
$\card{\aut G}=2^\w$.

Next note that $U_p\iso\Z\cross C_2$ for each prime $p$. Thus $\aut
G\iso P\cross V$, where $P$ is the direct product of countably many
copies of $\Z$ and $V$ is the direct product of countably many copies
of $C_2$. Since each nonzero element of $V$ has order $2$, it follows
that $V$ is isomorphic to a direct sum of $\card{V}$ copies of $C_2$.
Thus we can identify $V$ with a vector space of dimension $2^\w$ over
the field of two elements. Hence, $$\aut(\aut G)\iso \aut P\cross\aut
V=\aut P\cross GL(V),$$ where $GL(V)$ is the general linear group on
the vector space $V$. Since $\card{GL(V)}=2^{2^\w}$, it follows that
$\card{\aut(\aut G)}=2^{2^\w}$.\qed{ }

Enriqueta Rodr\'\i guez-Carrington has observed that the natural
modification of my argument shows that the derivation tower of every
Lie algebra eventually leads to a centerless Lie algebra. Since Simon
Thomas [1985] proved that the derivation tower of every centerless Lie
algebra must eventually terminate, it follows that the derivation tower
of any Lie algebra must eventually terminate.

One might hope, since every step of the automorphism tower kills the
center of the previous group, that $G_\omega$ is always centerless; but
this is not so. The dihedral group with eight elements has a center of
size two, but is isomorphic to its own automorphism group (there is an
outer automorphism which swaps $a$ and $b$ in the presentation
$\<a,b\st a^2=1,b^2=1,(ab)^4=1>$ ).  The group at stage $\omega$ is
just the two element group, which still has a center, and so this tower
survives until $\omega+1$. Simon Thomas has constructed examples
showing that for every natural number $n$ there are finite groups whose
tower has height $\omega+n$, but these also become centerless at stage
$\w+1$. Perhaps our attention should focus, for an arbitrary group $G$,
on the least ordinal stage $\g$ such that $G_\g$ is centerless. The
main point, then, is to find an explicit bound on how large $\g$ can be
in comparison with $\card{G}$.

To my knowledge, the following questions remain open: For which ordinals 
$\g$ is there a group whose tower becomes centerless in exactly
$\g$ many steps? Is there a countable group with an uncountable
automorphism tower? Is there a finite group with an uncountable
automorphism tower? Is there a finite group $G$ such that $G_\w$ is
infinite?

\bigskip
{\tenpoint\tightlineskip\obeylines
{\sc Joel David Hamkins, Mathematics
City University of New York, College of Staten Island 
Staten Island, NY 10314} 
{\tt hamkins@integral.math.csi.cuny.edu, www.library.csi.cuny.edu/users/hamkins/}} 

\bigskip
\nopagenumbers
\parindent=0pt
\newbox\Article
\newbox\Journal
\newbox\Author
\newbox\Vol
\newbox\No
\newbox\Year
\newbox\Page
\newbox\Book
\newbox\Publisher
\newbox\Pubaddr
\newbox\Key
\newbox\Editor
\newbox\Comment
\newbox\Note
\def\entry#1#2\par{\item{#1\quad}\hskip-1.1em#2\par}
\def\article#1{\setbox\Article=\hbox{\sl #1, }}
\def\journal#1{\setbox\Journal=\hbox{\rm #1 }}
\def\author#1{\setbox\Author=\hbox{\sc #1, }}
\def\vol#1{\setbox\Vol=\hbox{\bf #1 }}
\def\no#1{\setbox\No=\hbox{no. #1 }}
\def\year#1{\setbox\Year=\hbox{\rm({\oldstyle #1}) }}
\def\page#1{\setbox\Page=\hbox{\rm p. #1 }}
\def\book#1{\setbox\Book=\hbox{\it #1, }}
\def\publisher#1{\setbox\Publisher=\hbox{\rm #1, }}
\def\pubaddr#1{\setbox\Pubaddr=\hbox{\rm #1, }}
\def\key#1{\setbox\Key=\hbox{#1}}
\def\editor#1{\setbox\Editor=\hbox{\rm(#1, Ed.) }}
\def\comment#1{\setbox\Comment=\hbox{\rm #1}}
\def\note#1{\setbox\Note=\hbox{\rm #1 }}
\def\ref#1\par{\smallskip{#1
  \entry{\ifhbox\Key\unhbox\Key\else[\ ]\fi}%
  \unhbox\Author\unhbox\Note
  \ifhbox\Book \unhbox\Book\unhbox\Publisher\unhbox\Pubaddr
               \unhbox\Editor\unhbox\Page\unhbox\Year\unhbox\Comment
  \else \unhbox\Article\unhbox\Journal\unhbox\Vol\unhbox\No\unhbox\Editor
        \unhbox\Page\unhbox\Year\unhbox\Comment\fi\par}}

\tenpoint

\ref
\author{Joel David Hamkins and Simon Thomas}
\article{Changing the heights of automorphism towers}
\journal{submitted to the Annals of Pure and Applied Logic}
\key{[1997]}

\ref
\author{J. A. Hulse}
\article{Automorphism towers of polycyclic groups}
\journal{Journal of Algebra}
\vol{16}
\year{1970}
\page{347--398}
\key{[1970]}

\ref
\author{Andrew Rae and James E. Roseblade}
\article{Automorphism Towers of Extremal Groups}
\journal{Math. Z.}
\page{70-75}
\year{1970}
\vol{117}
\key{[1970]}

\ref
\author{Simon Thomas}
\article{The automorphism tower problem}
\journal{Proceedings of the American Mathematical Society}
\vol{95}
\year{1985}
\page{166--168}
\key{[1985]}

\ref
\author{Simon Thomas}
\article{The automorphism tower problem II}
\journal{to appear in Israel J. Math.}
\key{[1997]}

\ref
\author{H. Wielandt}
\article{Eine Verallgemeinerung der invarianten Untergruppen}
\journal{Math. Z.}
\vol{45}
\year{1939}
\page{209--244}
\key{[1939]}

\bye